\def\simplelatex{\iftrue}
\documentclass[12pt]{article} 
\usepackage[francais]{babel}
\simplelatex%
\else%
\usepackage[latin1]{inputenc}
\usepackage[T1]{fontenc}
\fi%
\usepackage{amsmath}
\usepackage{amsfonts}
\usepackage{amssymb}
\usepackage{amsthm}
\simplelatex%
\let\mathscr=\mathcal
\else%
\usepackage{mathrsfs}
\fi%
\simplelatex%
\let\url=\texttt
\else%
\usepackage{url}
\fi%
\newenvironment{rmq}[1][]{\refstepcounter{prop}%
\bigskip
\noindent{\textbf{Remarque \theprop{}. #1}}}{}

\newtheorem{prop}{Proposition}[section]
\newtheorem{thm}[prop]{Th\'eor\`eme}

\newcommand{\br}{\mathop\mathrm{Br}\nolimits}

\newcommand{\bP}{{\mathbb P}}
\newcommand{\Q}{{\mathbb Q}}

\newcommand{\Z}{{\mathbb Z}}

\simplelatex%
\let\bimu=\mu
\let\mathbi=\mathbf
\else%
\DeclareFontFamily{OML}{cmmib}{\skewchar\font127 }
\DeclareFontShape{OML}{cmmib}{m}{it}%
       {<5><6><7><8><9>gen*cmmib%
        <10><10.95>cmmib10%
        <12><14.4><17.28><20.74><24.88>cmmib12%
        }{}
\DeclareSymbolFont{biletters}{OML}{cmmib}{m}{it}
\DeclareSymbolFontAlphabet{\mathbi}{biletters}
\DeclareMathSymbol{\bimu}{\mathord}{biletters}{"16}
\fi%
\simplelatex%

\fi%

\begin{document}
\title{Remarques sur un article r\'ecent de B. Poonen}
\author{ J-L. Colliot-Th\'el\`ene}
\maketitle 

  {\em R\'esum\'e} :    B. Poonen a r\'ecemment exhib\'e des exemples de vari\'et\'es projectives  
 et lisses  de dimension 3 sur un corps de nombres qui n'ont pas
  de point rationnel et pour lesquelles il n'y a pas
 d'obstruction de Brauer--Manin apr\`es rev\^etement fini \'etale.
 Dans cette note, je montre que les vari\'et\'es qu'il construit poss\`edent
 des z\'ero-cycles de degr\'e 1. 
 
 \smallskip
 
  {\em Summary} : B. Poonen recently produced smooth threefolds
  over a number field which do not have a rational point
 but have no Brauer--Manin obstruction even after descent to
 a finite \'etale cover. In this note I show that the varieties he produces
 have   zero-cycles of degree 1.
 
 \bigskip
 
  {\em Mots-cl\'es} : points rationnels, z\'ero-cycles, principe de Hasse, obstruction de Brauer--Manin  


Mathematical Subject Classification (MSC 2000) 
 Primary :  
14G05, 14G25 11G35. 
 Secondary : 14J20, 14F22.

 \bigskip
 
 \section{Introduction}

Soit $k$ un corps de nombres.  Dans son article \cite{Poonen},   B. Poonen construit des exemples de vari\'et\'es projectives  
 et lisses    de dimension 3 sur $k$  qui ont les propri\'et\'es suivantes :
 elles ont des points dans tous les compl\'et\'es de $k$, il n'y a pas d'obstruction de Brauer--Manin
 \`a l'existence d'un point rationnel, mieux, il n'y a pas d'obstruction de Brauer-Manin apr\`es
 descente par rev\^etements finis \'etales,
  et cependant les vari\'et\'es  ne poss\`edent  pas de point $k$-rationnel. 
  C. Demarche \cite{Demarche} 
  vient de montrer que ceci implique qu'aucune obstruction de descente sous un groupe alg\'ebrique ne saurait  rendre compte de l'inexistence de point rationnel. La situation est donc tout \`a fait diff\'erente
 de celle de l'exemple  historique de Skorobogatov (\cite {Sko}, \cite{HaraSko}).

 Dans \cite{CTconj} j'ai conjectur\'e que pour toute vari\'et\'e projective, lisse, g\'eom\'etri\-quement connexe sur un corps
 de nombres, l'obstruction de Brauer-Manin \`a l'existence d'un z\'ero-cycle de degr\'e 1 est la seule obstruction.
 
 Dans cette note je montre que la conjecture vaut pour les vari\'et\'es cons\-truites par Poonen : 
 elles ont toutes  un z\'ero-cycle de degr\'e 1.  La question reste enti\`ere pour l'exemple de Skorobogatov.
 
 Au paragraphe 1 je rappelle bri\`evement la construction des vari\'et\'es.
 Au paragraphe  2 je donne une d\'emonstration alternative du calcul  
 du groupe de Brauer de ces vari\'et\'es. Le r\'esultat de ce calcul permet de montrer (\cite{Poonen})
 qu'il n'y a pas d'obstruction de Brauer--Manin  \`a l'existence d'un point rationnel,
 et {\it a fortiori} pas d'obstruction de Brauer--Manin \`a l'existence de z\'ero-cycles de degr\'e 1.
 Au paragraphe 3, logiquement
 ind\'ependant du pr\'ec\'edent, j'\'etablis l'existence de z\'ero-cycles de degr\'e 1.

 Cet article a \'et\'e con\c cu  \`a l'occasion d'expos\'es    donn\'es en avril 2008 \`a l'universit\'e Emory (Atlanta, G\'eorgie).

\section{Les vari\'et\'es}

Soit $k$ un corps de caract\'eristique nulle.
Les vari\'et\'es consid\'er\'ees dans  \cite{Poonen}   s'ins\`erent dans le diagramme suivant, que nous allons d\'etailler.
 
$$\begin{array}{ccccccccccccccccccccccccccccccccccccccccccccccccccccc}
 X & \to &  V\cr
\downarrow & & \downarrow \cr
 C \times \bP^1 & \to & \bP^1 \times \bP^1 \cr
\downarrow & &  \downarrow \cr
 C&  \to &  \bP^1
 \end{array} $$

Dans ce diagramme, toutes les vari\'et\'es sont des $k$-vari\'et\'es projectives, lisses, g\'eom\'etriquement connexes.
La vari\'et\'e $C$ est une courbe.
Les fl\`eches verticales inf\'erieures  sont donn\'ees par la projection
sur le premier facteur.
La vari\'et\'e $V$ est de dimension 3, elle est fibr\'ee en coniques sur $\bP^1 \times \bP^1$.
Le lieu de ramification de cette fibration (le lieu o\`u la conique fibre est singuli\`ere)
est une $k$-courbe projective, lisse, g\'eom\'etriquement int\`egre $Z_{1} \subset \bP^1 \times \bP^1$.
La projection sur le premier facteur $ \bP^1 \times \bP^1 \to \bP^1$ induit un rev\^etement
ramifi\'e $Z_{1} \to \bP^1$, le lieu de ramification dans $\bP^1$ ne contient pas le
point $\infty \in \bP^1(k)$.

Le morphisme $C \to \bP^1$ est un rev\^etement ramifi\'e, son lieu
de ramification dans $\bP^1$ ne rencontre pas le lieu de ramification
de $Z_{1} \to \bP^1$.
 
La partie verticale gauche du diagramme est obtenue par r\'etrotirette 
\`a partir de la partie verticale droite (produits fibr\'es 
via le morphisme $C \to \bP^1$).

Pour tout point sch\'ematique $w\in \bP^1$,  de corps r\'esiduel $k(w)$, la fibre $V_{w}$ est une surface 
g\'eom\'e\-tri\-quement int\`egre sur $k(w)$ qui contient un ouvert affine d'\'equation  
$$y^2-az^2=P_{w}(x)$$
avec $a \in k \setminus k^2$ et $P_{w}(x) \in k(w)[x]$ un polyn\^ome non nul, de degr\'e 4 en la variable $x$,
s\'eparable pour $w$ non dans le lieu de ramification de $Z_{1} \to \bP^1$.

Le lieu de d\'eg\'en\'erescence de la fibration en coniques $X \to C \times \bP^1$
est l'image inverse $Z$ de $Z_{1}$ par $C \times \bP^1 \to \bP^1 \times \bP^1$.
Il est montr\'e dans \cite{Poonen}  (Lemma 8.1) que $Z \subset C \times \bP^1$
est une $k$-courbe projective, lisse, g\'eom\'etriquement int\`egre.

La description des fibres de l'application compos\'ee $V \to \bP^1 \times \bP^1 \to \bP^1$
(premi\`ere projection)
  montre que l'on a :

\begin{prop}\label{fibresgeomint}
Sous les hypoth\`eses faites ci-dessus, pour tout point sch\'e\-ma\-tique $P$
de $C$, la fibre $X_{P}/k(P)$ du morphisme   $X    \to C$ au point $P$
est une  $k(P)$-vari\'et\'e g\'eom\'etriquement int\`egre.
\end{prop}

\section{Leur groupe de Brauer}

Le r\'esultat suivant est \'etabli dans \cite{Poonen} au moyen d'une \'etude de l'action du groupe de Galois sur le 
groupe de Picard g\'eom\'etrique de $X$. Nous proposons une d\'emonstration un peu plus courte.
La diff\'erence entre la d\'emonstration de \cite{Poonen} et la pr\'esente d\'emonstration est essentiellement la m\^eme
que celle entre les propositions 7.1.1  et 7.1.2 de \cite{Skolivre}.

\begin{prop}\label{brauer}
Sous les hypoth\`eses ci-dessus, la fl\`eche naturelle de groupes de Brauer
$\br C \to \br X$ est un isomorphisme.
\end{prop}

\begin{proof}

Notons $B=C \times \bP^1$. La fibration en coniques $X \to B$ est d\'eg\'en\'er\'ee en un seul point de
codimension 1 de $B$, correspondant \`a la courbe g\'eom\'etriquement int\`egre $Z$.

Soit $\eta$ le point g\'en\'erique de $B$ et $X_{\eta}$ la fibre g\'en\'erique de $X \to B$.
C'est une conique lisse.

De fa\c con g\'en\'erale, \'etant donn\'e un $k$-morphisme plat $X  \to B$
de $k$-vari\'et\'es lisses  g\'eom\'etriquement int\`egres,
on dispose d'un diagramme commutatif de suites exactes de groupes
de cohomologie \'etale :

$$\begin{array}{ccccccccccccccccccccccccccccccccccccccccccccccccccccc}
0 & \to & \br B  & \to &  \br k(B)  &\buildrel{  \oplus_{x}   {\partial_x}     } \over \longrightarrow  &  \bigoplus_{x  \in B^{(1)}}  H^1(k(x),{\Q}/{\Z}) \cr
&   & \downarrow &&\downarrow{}{    }     & & \downarrow{}{e_{y/x}.{\rm Res}_{k(x)/k(y)}} \cr
0 & \to & \br  X & \to & \br X_{\eta}  &\buildrel{  \oplus_{x}   {\partial_y}     }  \over \longrightarrow   &      \bigoplus_{x\in B^{(1)}   } \bigoplus_{y\in X^{(1)}, y \to x }    H^1(k(y) ,\Q/\Z).     \cr
\end{array} $$ 
Dans ce diagramme, $x$ parcourt les points de codimension 1 de $B$, et $y$ parcourt les points de codimension 1
de $X$ qui ne sont pas situ\'es sur $X_{\eta}$. Pour $y \in X$ d'image $x \in B$, on a l'inclusion
des corps r\'esiduels $k(x) \subset k(y)$. L'entier $e_{y/x}$ est l'indice de ramification.
Les fl\`eches $\partial$ sont les fl\`eches de r\'esidu.

On peut extraire un tel diagramme des expos\'es de Grothendieck sur le groupe de Brauer 
(\cite{Grothendieck}, voir GB II, Cor. 1.10,  GB III, Prop. 2.1,  GB III, Thm. 6.1).
Pour plus de d\'etails, voir  par exemple \cite[\S 3]{CTBarbara} et \cite[Chapter 6]{GilleSzamuely}. 
La d\'emonstration combine les suites de localisation, leur fonctorialit\'e et
le th\'eor\`eme de puret\'e pour le groupe de Brauer.

Dans la pr\'esente situation, la fibre g\'en\'erique est une conique sur $k(B)$. Soit $\beta \in \br k(B)$
la classe de l'alg\`ebre de quaternions associ\'ee \`a cette conique. 
Cette classe $\beta$ est non triviale, elle admet un unique r\'esidu non trivial, au point g\'en\'erique de $Z$,
et la classe correspondante est la classe de $a \in k^*/k^{*2} = H^1(k,\Z/2) \subset H^1(k(Z),\Q/\Z)$.

Pour la conique $X_{\eta}$ sur $k(B)$ sans point rationnel, on dispose de la suite exacte
classique  (cf. \cite{Sko} p. 138)
$$0 \to \Z/2.{\beta} \to \br k(B) \to \br X_{\eta} \to 0.$$
Le diagramme ci-dessus montre alors que l'application $\br B \to \br X$ est injective.

Soit $\alpha \in \br X$.  L'image de $\alpha$ dans $\br X_{\eta}$ est l'image d'un \'el\'ement
$\gamma  \in \br k(B)$,   d\'efini \`a addition pr\`es de $\beta$.

Le lieu de ramification de $X \to B$ est r\'eduit \`a l'unique courbe $Z \subset B$.
Ceci implique que pour $x \in B^{(1)}$ diff\'erent du point g\'en\'erique de $Z$, et $y$ l'unique point
de $X^{(1)}$ au-dessus de $x$, qui d\'efinit une conique lisse donc g\'eom\'etriquement int\`egre sur
le corps $k(x)$,
la fl\`eche $e_{y/x}.{\rm Res}_{k(x)/k(y)}={\rm Res}_{k(x)/k(y)} : H^1(k(x), \Q/\Z) \to H^1(k(y), \Q/\Z)$ est une injection.
Quant  \`a la fl\`eche $e_{y/x}.{\rm Res}_{k(x)/k(y)}= {\rm Res}_{k(x)/k(y)}$ 
associ\'ee au point g\'en\'erique de $Z$ et \`a l'unique $y$ au-dessus de ce point, 
 son noyau est le noyau de la fl\`eche de restriction
$H^1(k(Z), \Q/\Z) \to H^1(k(Z)(\sqrt{a}), \Q/\Z)$. Ce noyau est d'ordre 2, engendr\'e par
la classe de $a$ dans $ k^*/k^{*2} =H^1(k,\Z/2)$.

Du diagramme ci-dessus on conclut que les r\'esidus de $\gamma$
aux points autres que le point g\'en\'erique de $Z$ sont  nuls, et qu'au point
g\'en\'erique de $Z$ le r\'esidu est soit   trivial soit    \'egal au r\'esidu de $\beta$.
Quitte \`a remplacer \'eventuellement $\gamma$ par $\gamma+\beta$, ce
qui est loisible, on voit que $\gamma \in \br k(B)$ est dans l'image de $\br B$.
Comme  l'application $\br X  \to \br X_{\eta}$ est injective, ceci ach\`eve de
montrer que l'application $\br B \to \br X$ est un isomorphisme.

Le groupe de Brauer de la droite projective sur un corps est \'egal
\`a l'image du groupe de Brauer de ce corps. De $B = C \times \bP^1$
on d\'eduit que $\br C \to \br B$ est un isomorphisme. Ceci ach\`eve la
d\'emonstration de la proposition.
\end{proof}

\section{Existence de z\'ero-cycles de degr\'e 1}

\begin{thm}
Soit $k$ un corps de nombres.
Soient $V \to \bP^1 \times \bP^1$, $C \to \bP^1$ et $X \to C \times \bP^1$ comme 
ci-dessus. Faisons les deux hypoth\`eses :

(i)  la fibre (lisse)  $V_{\infty}$ de la fibration compos\'ee
$V \to  \bP^1 \times \bP^1 \to  \bP^1$ (premi\`ere projection)
a des points dans tous les compl\'et\'es de $k$;

(ii) il existe un $k$-point $P$ de $C$ dont l'image par $C \to \bP^1$
est le point $\infty$.

\noindent Alors pour tout entier
  $n \geq 2g+1$, o\`u $g$ d\'esigne le genre de $C$,
  la $k$-vari\'et\'e $X$ poss\`ede un point dans un corps
  extension de $k$ de degr\'e $n$. En particulier  $X$ 
  poss\`ede un z\'ero-cycle de degr\'e 1. 
\end{thm}

 \begin{proof}
 Soit $n$ comme dans l'\'enonc\'e. Le syst\`eme lin\'eaire associ\'e au diviseur $nP$
 de $C$ est alors tr\`es ample. Soit $C \hookrightarrow {\bP}^N$ le plongement associ\'e,
 et soit  $\check{\bP}^N$ l'espace projectif dual de $\bP^N$. 
 
 La vari\'et\'e d'incidence $W \subset  \bP^N\times \check{\bP}^N$ d\'efinie par l'annulation de la  forme
 bihomog\`ene $\sum_{i=0}^NX_{i}Y_{i}$
 est un diviseur tr\`es ample. La projection $W \to  \check{\bP}^N$ d\'efinit un
 fibr\'e (lisse) en espaces projectifs $\bP^{N-1}$.
 
   Soit $W_{Z} \subset Z \times \check{\bP}^N$ l'image inverse de la correspondance d'incidence
 via le morphisme compos\'e 
 $$Z \times \check{\bP}^N \to  C \times  \check{\bP}^N  \to \bP^N  \times \check{\bP}^N .$$ 
 C'est un fibr\'e (lisse) en  espaces projectifs $\bP^{N-1}$ au-dessus
 de la courbe $Z$, qui est projective, lisse et g\'eom\'etriquement int\`egre. 
 Ainsi $W_{Z}$ est projectif, lisse et g\'eom\'etriquement int\`egre.
 
 Le m\^eme argument montre que l'image inverse $W_{C}\subset C \times \check{\bP}^N$
 de la correspondance d'incidence  est une $k$-vari\'et\'e projective, lisse, g\'eom\'etriquement int\`egre.
 
 Les morphismes naturels $W_{Z} \to W_{C}  \to \check{\bP}^N$ sont des morphismes 
 propres dominants et finis -- ils sont quasi-finis  car la courbe
  $C$ engendre $\bP^N$ projectivement.
  
 \medskip
  
  La proposition \ref{fibresgeomint} et un \'enonc\'e bien
 connu (\cite[I, p.~43]{CTSaSwD},  \cite[Thm. 1]{SkoSTNparis}, \cite[\S 3.2]{CTToulouse}) 
  impliquent l'existence d'un ensemble fini $S$ de places
 de $k$ tel que pour toute extension finie $L$ de $k$
 et toute place $w$ de $L$ non situ\'ee au-dessus d'une place de $S$,
 l'application $X(L_{w}) \to C(L_{w})$ est surjective.

\medskip
  
  Pour toute place $v$ de $k$,  la fibre (lisse) $X_{P}$ du morphisme $X \to C$ 
  au-dessus du $k$-point $P$ poss\`ede un $k_{v}$-point lisse.
  Ceci implique que pour toute extension finie $L/k_{v}$,  l'image de
  l'application induite $X(L) \to C(L)$ contient un 
  voisinage ouvert de $P\times_{k_{v}}L$
  pour la topologie
$v$-adique sur $C(L)$.

 \medskip
 
 Il existe un hyperplan $H_{0} \in  \check{\bP}^N(k)$ qui  d\'ecoupe exactement
 le diviseur $nP$ sur $C$. Soit $v$ une place de $k$.
 Tout hyperplan $H \in  \check{\bP}^N(k)$
 qui est suffisamment proche de $H_{0}$ pour la topologie $v$-adique
 d\'ecoupe sur la courbe $C$ une somme de points ferm\'es
 qui lorsque l'on passe de $k$ \`a un compl\'et\'e $k_{v}$ donnent des points
 sur diverses extensions finies $L$ de $k_{v}$  de degr\'e au plus $n$,
points qui sont proches de $P\times_{k}L$ pour la topologie
$v$-adique sur $C(L)$. Rappelons qu'un tel corps local $k_{v}$ (de caract\'eristique
nulle) n'admet   qu'un nombre fini
d'extensions   $L/k_{v}$ de degr\'e au plus    $n$.

\medskip

Le th\'eor\`eme d'irr\'eductibilit\'e de Hilbert avec approximation faible
(cf. \cite{Harari}, Prop. 3.2.1), appliqu\'e \`a l'ensemble fini $S$
de places et 
au rev\^etement fini de vari\'et\'es
g\'eom\'etriquement int\`egres $W_{Z} \to  \check{\bP}^N$,
compos\'e de $W_{Z} \to  W_{C} \to  \check{\bP}^N$,
permet alors de trouver un point ferm\'e $M \in C$,
de degr\'e $n$ sur $k$, tel que la fibre $X_{M}/k(M)$ soit  une surface de
Ch\^atelet lisse, poss\`ede
des points dans tous les compl\'et\'es de $k(M)$
(ceci par combinaison des trois paragraphes pr\'ec\'edents)
et que de plus le sch\'ema $Z_{M} = Z \times_{C}M$ 
soit  int\`egre. Mais ceci signifie que la surface de Ch\^atelet
correspondante admet une \'equation affine
$y^2-az^2=P(x)$ avec $P(x)\in k(M)[x]$ irr\'eductible de degr\'e 4.
L'un des principaux r\'esultats de  \cite{CTSaSwD}, le th\'eor\`eme
  8.11,  assure que
le principe de Hasse vaut pour une telle surface. Ainsi $X_{M}(k(M)) \neq
\emptyset$, {\it a fortiori} $X(k(M)) \neq \emptyset$.

\end{proof}

 \begin{rmq}
Supposons $C=\bP^1$. Dans ce cas on peut prendre $n=1$, et on obtient
  $X(k)\neq \emptyset$.  On a mieux. Comme les surfaces de Ch\^atelet d'\'equation affine
 $y^2-az^2=P(x)$ avec $P(x)$ irr\'eductible satisfont l'approximation
 faible (\cite[Thm. 8.11]{CTSaSwD}), une variante imm\'ediate de l'argument  \'etablit     l'approximation faible pour $X$: 
 l'ensemble $X(k)$ est dense dans le produit des $X(k_{v})$.
 L'argument est   un cas particulier d'un th\'eor\`eme g\'en\'eral
 (Harari,  \cite{Harari}, Thm. 4.2.1).
      \end{rmq}

 \begin{rmq} Dans \cite{CTconj}, outre la conjecture sur l'existence de z\'ero-cycles de degr\'e 1,
  j'ai   propos\'e une conjecture   
 sur  l'application diagonale de groupes de Chow de 
 z\'ero-cycles  
 $ CH_{0}(X) \to \prod_{v} CH_{0}(X_{v}).$
 Dans une s\'erie d'articles (voir la bibliographie de \cite{CTconj}), des versions de cette conjecture ont \'et\'e 
 \'etablies  pour des vari\'et\'es fibr\'ees en vari\'et\'es de Severi-Brauer
 au-dessus d'une courbe. Il serait int\'eressant de voir si l'on peut
 pousser les arguments ci-dessus et \'etablir la conjecture pour les familles $X \to C$
 de surfaces de Ch\^atelet consid\'er\'ees ici.
  \end{rmq}

\bigskip
 Jean-Louis Colliot-Th\'el\`ene 
 
 C.N.R.S., UMR  8628,

 Math\'ematiques,
 B\^atiment 425, 
 
 Universit\'e   Paris-Sud,
 
F-91405 Orsay, 

France 

\medskip

courriel :  jlct \`a math.u-psud.fr 
\end{document}